\documentclass[12pt]{amsart}

\usepackage[margin=1in]{geometry}

\usepackage{amsmath,amscd,amssymb,amsfonts,latexsym,wasysym, mathrsfs, mathtools,hhline,color}
\usepackage[all, cmtip]{xy}
\usepackage{csquotes}
\usepackage{url}
\usepackage{comment}
\usepackage{tikz}
\usepackage{caption}

\definecolor{hot}{RGB}{65,105,225}

\usepackage[pagebackref=true,colorlinks=true, linkcolor=hot ,  citecolor=hot, urlcolor=hot]{hyperref}%make all the references and links clickable
\usepackage{ textcomp }
\usepackage{ tipa }
\usepackage{graphicx,enumerate}

\usepackage{gensymb}
\usepackage{listings,xcolor}
\usepackage{geometry}
\theoremstyle{plain}
\newtheorem{theorem}{Theorem}[section]
\newtheorem{prop}[theorem]{Proposition}

\newtheorem{lm}[theorem]{Lemma}

\newtheorem{cor}[theorem]{Corollary}

\newtheorem{thrm}[theorem]{Theorem}

\theoremstyle{definition}

\newtheorem{defn}[theorem]{Definition}
\newtheorem{que}[theorem]{Question}
\newtheorem{rmk}[theorem]{Remark}

\newtheorem{ex}[theorem]{Example}
\newtheorem*{ex*}{Example}

\def\be{\begin{equation}}
\def\ee{\end{equation}}

\def\bt{\begin{thrm}}
\def\et{\end{thrm}}

\def\bc{\begin{cor}}
\def\ec{\end{cor}}

\def\br{\begin{rmk}}
\def\er{\end{rmk}}

\def\bp{\begin{prop}}
\def\ep{\end{prop}}

\def\bl{\begin{lm}}
\def\el{\end{lm}}

\def\bex{\begin{ex}}
\def\eex{\end{ex}}

\def\bd{\begin{defn}}
\def\ed{\end{defn}}

\newcommand{\C}{\mathbb{C}}
\newcommand{\Z}{\mathbb{Z}}
\newcommand{\Hom}{\mathrm{Hom}}

\newcommand{\Cha}{\mathrm{Char}}

\newcommand{\orb}{\mathrm{orb}}

\newcommand{\homo}{\mathrm{Hom}}
\newcommand{\sV}{\mathcal{V}}

\newcommand{\sA}{\mathcal{A}}

\newcommand{\sB}{\mathcal{B}}

\newcommand{\CP}{\mathbb{CP}}

\newcommand{\bm}[1]{\mbox{\boldmath{$#1$}}}

\title[Double star arrangement and the pointed multinet]{Double star arrangement and the pointed multinet}

\author{Yongqiang Liu}
\address{Yongqiang Liu: Institute of Geometry and Physics, University of Science and Technology of China, 96 Jinzhai Road, Hefei Anhui 230026 China} 
\email{liuyq@ustc.edu.cn}

\author{Wentao Xie}
\address{Wentao Xie: School of Mathematical Sciences, University of Science and Technology of China, 96 Jinzhai Road, Hefei, Anhui 230026 China}
\email{xwt@mail.ustc.edu.cn}

\date{\today}

\keywords{Hyperplane arrangement, cohomology jump loci, orbifold fibration, pointed multinet}

\subjclass[2010]{14F45, 52C35, 32S55}

\begin{document}

\begin{abstract}
Let $\sA$ be a hyperplane arrangement in a complex projective space. It is an open question if the degree one cohomology jump loci (with complex coefficients) are determined by the combinatorics of $\sA$. By the work of Falk and Yuzvinsky \cite{FY}  and Marco-Buzun\'ariz \cite{Mar},  all the irreducible components passing through the origin are determined by the multinet structure, which is combinatorially determined.  Denham and Suciu introduced the pointed multinet structure, which is combinatorially determined, to obtain examples of arrangements with  translated positive-dimensional components in the degree one cohomology jump loci \cite{DS}.
Suciu asked the question if all translated positive-dimensional components appear in this manner \cite{Suc14}. In this paper, we show that  the double star arrangement introduced by  Ishibashi, Sugawara and Yoshinaga \cite[Example 3.2]{ISY22}  gives a negative answer to this question. 
\end{abstract}

\maketitle

\section{Background}
Let $\sA$ be an arrangement of finite hyperplanes in $\CP^d$ with $d\geq 2$. 
The combinatorics of the arrangement $\sA$ is encoded in its intersection lattice, $L(\sA)$; this is the poset of all intersections of $\sA$, ordered by reverse inclusion, and ranked by codimension. Let $X(\sA)$ denote the complement of the hyperplane arrangement $\sA$. A fundamental problem in the theory of hyperplane arrangements is to decide whether various invariants of $X(\sA)$ are determined by  the combinatorics $L(\sA)$.  For instance, Betti numbers and the cohomology ring of $X(\sA)$ are combinatorially determined (e.g., see \cite{OS}), while the fundamental group of $X(\sA)$ cannot be determined by $L(\sA)$ due to Rybnikov \cite{Ryb}.  See \cite{Dim17} for an overview of the theory.
%It is an open question whether cohomology jump loci of hyperplane arrangement complement are combinatorially determined; see \cite{Suc14}  for recent progress in this direction, and also \cite{Dim17} for an overview of the theory. 

\subsection{Cohomology jump loci}
In this paper, we focus on the degree one cohomology jump loci. We first recall its definition.

Assume that $X$ is a connected finite CW complex.
Let $\Cha(X)$ denote the group of $\C^*$-valued  characters $ \homo(\pi_1(X),\C^*)$, which is a commutative affine algebraic group. Each character $\rho \in \Cha(X)$ defines a rank one $\C$-local system on $X$, denoted by $L_{\rho}$. 
\bd
The degree one cohomology jump loci of $X$ is defined as
$$\sV^1(X)\coloneqq \lbrace \rho\in \Cha(X) \mid  H^{1}(X, L_{\rho})\neq 0 \rbrace.$$
\ed
Since $\sV^1(X)$ only depends on the the fundamental group of $X$, $\sV^1(G)$ is well-defined for any finitely presented group $G$.

For hyperplane arrangement complement $X(\sA)$, we simply write $\sV^1(\sA)$ for $\sV^1(X(\sA))$.
Since we are only interested in $\sV^1(\sA)$, 
using Lefschetz hyperplane section theorem, {\bf we always assume that 
$\sA$  is  an essential  line arrangement in $\mathbb{CP}^2$.} 
Being essential means that the intersection of all lines of $\sA$ is empty. 
The following question is widely open, see e.g. \cite[Problem 3.15]{Suc14}.  
\begin{center}
Is $\sV^1(\sA)$ determined by $L(\sA)$?
\end{center}

\medskip

We now describe in more details about $\sV^1(\sA)$.  A  celebrity result  due to  Arapura \cite{Ara97} and Artal Bartolo, Cogolludo-Agust\'in and Matei \cite{ACM13} shows that $\sV^1(X)$  are unions of torsion-translated
subtori. To explain their results, we recall the notion of orbifold maps.

 An algebraic map $f\colon X(\sA) \to \Sigma_k$  is called an (genus 0) orbifold fibration, if $f$ is surjective, has connected generic fiber and $\Sigma_k=\CP^1 -\{k \text{ points}\}$ with $k\geq 2$. 
There exists a maximal Zariski open subset $U\subset \Sigma_k$ such that $f$ is a fibration over $U$. Say $B=\Sigma_k-U$ (could be empty) has $s$ points, denoted by $\{q_1,\ldots,q_s\}$.  We assign
the multiplicity $\mu_j$ of the fiber   $f^{*}(q_j)$  (the $\gcd$ of the coefficients of the divisor $f^* q_j$) to the point $q_j$. Set  ${\bm \mu}=(\mu_1,\ldots,\mu_s)$. We only count such multiplicity if $\mu_j\geq 2$.  %The orbifold $\Sigma$ is said to be maximal if no divisor $f^{\ast}(p)$ is an $n$-multiple for $n>\mu(p)$.
%Such orbifold map $f$ is called of type $(k,{\bm \mu})$, where ${\bm \mu}=(\mu_1,\ldots,\mu_s)$. 

The orbifold group $\pi_1^{\orb}(\Sigma_k,{\bm \mu})$ associated to these data is  defined as
\[ \pi_{1}^{\orb}(\Sigma_k,{\bm \mu})\coloneqq F_{k-1} \ast\Z_{\mu_1}\ast \cdots \ast\Z_{\mu_s}, \]
where $F_{k-1}$ is the free group with $(k-1)$ generators. 
The orbifold fibration  $f$  induces a surjective map to the orbifold group 
$ \pi_1(X(\sA)) \twoheadrightarrow \pi_1^{\orb}(\Sigma_k,{\bm \mu}),$
which gives an embedding (see e.g. \cite[Proposition A.1]{Suc14})
$$f^*\colon \sV^1(\pi_{1}^{\orb}(\Sigma_k,{\bm \mu})) \to \sV^1(\sA) .$$ 

Set $\pi=\pi_{1}^{\orb}(\Sigma_k,{\bm \mu})$, $A= \Z_{\mu_1}\oplus \cdots \oplus\Z_{\mu_s}$ and $\Cha(A)=\Hom(A,\C^*)$. Note that  $\Cha(\pi)=\Cha(\pi)^\circ\times \Cha(A)$,  where $\Cha(\pi)^\circ\cong (\C^*)^{k-1}$ 
is the identity component of $\Cha(\pi)$. A Fox calculus computation shows that 
\begin{center}
$\sV^1(\pi)=\begin{cases}
\big(\Cha(\pi)\setminus \Cha(\pi)^\circ \big)\cup \{\bm 1\}, & \text{ if } k=2  \text{ and } f \text{ has at least one multiple fiber,}\\
 \Cha(\pi), & \text{ if } k\geq 3,
\end{cases}  $
\end{center}
where we call $f$ small in the first case, and large
in the second case. See \cite[Section 3.4]{Suc14} for more details.

Following \cite[Theorem 3.8]{Suc14}, we summarize the results as follows.
\bt \cite{Ara97,ACM13} \label{thm structure} Let $\sA$  be  an essential  line arrangement in $\mathbb{CP}^2$.
Then  we have 
$$ \sV^1(\sA)=\bigcup_{f \text{ large} } f^*( \sV^1(\pi_{1}^{\orb}(\Sigma_k,{\bm \mu})) )\cup \bigcup_{f \text{ small} } f^*( \sV^1(\pi_{1}^{\orb}(\Sigma_k,{\bm \mu})) ) \cup Z,$$
where the unions are over the equivalence classes of orbifold fibration $f \colon X(\sA)\to \Sigma_k$ of the
types indicated, 
and $Z$ is a finite union of isolated torsion points. In particular, $\sV^1(\sA)$ has translated positive-dimensional components if and only if at least one of the orbifold fibration $f$ appearing in the union has multiple fibers.
\et 
\br This theorem indeed holds for any 
 complex smooth quasi-projective variety, see   \cite{Ara97,ACM13} and  also Budur and Wang's work \cite{BW15} for higher degree cohomology jump loci.
\er

\subsection{Multinet}
 As shown by Falk and Yuzvinsky in \cite{FY}, the
irreducible components of $\sV^1(\sA)$ passing through the origin can be described in terms
of multinet, which only depends on the intersection lattice  $L(\sA)$. 
Let us recall its definition.
\bd \label{def multinet}   A multinet on a line arrangement  $\sA$ in $\CP^2$ is a partition of $\sA$ into $k\geq 3$ subsets $\sA_1,\cdots, \sA_k$, together with an assignment of multiplicities, $m\colon \sA \to \Z_{>0}$, and a subset $\mathfrak{X}$ of intersection points in $\sA$, called the base locus, satisfying the following conditions:
\begin{itemize}
\item[(a)] $\sum_{H\in \sA_i} m_H= \kappa$, independent of $i$;
\item[(b)] for each $H\in \sA_i$ and $H'\in \sA_j$ with $i\neq j$, $H\cap H' \in \mathfrak{X}$; 
\item[(c)] for each $x\in \mathfrak{X}$, the sum $n_x\coloneqq\sum_{H\in \sA_i, x\in H} m_H$ is independent of $i$;
\item[(d)] for every $1\leq i \leq k$ and $H,H'\in \sA_i$, there is a sequence $\{H_0,H_1,\cdots,H_s\}$ with $H_0=H$ and $H_s=H'$ such that $H_{j-1}\cap H_j \notin \mathfrak{X}$ for all $1\leq j \leq s$;
\item[(e)] $\gcd(m_H\colon H\in \sA)$=1.
\end{itemize} Such multinet is called a $(k,\kappa)$-multinet.
\ed

The notion of multinet was simultaneously introduced by  Marco-Buzun\'ariz in \cite{Mar} by the name of combinatorial pencil. %Both Falk and Yuzvinsiky \cite{FY} and  Marco-Buzun\'ariz  discovered the relations between resonance varieties of the line arrangements and multinet. However we will not pursue this direction in this paper. 
 Next we list some known results about multinet due to Pereira and Yuzvinsky  as follows.
\bt \cite{PY,Yuz}\label{thm 34} Let $(\sA_1,\cdots,\sA_k)$ be a multinet on $\sA$ with a base locus $\mathfrak{X}$. 
\begin{itemize}
\item[(1)] If $\vert\mathfrak{X} \vert>1 $, then $k = 3$ or $4$.
\item[(2)] If there exists $H\in \sA$ such that $m_H>  1$, then $k = 3$.
\end{itemize}
\et

Let $f_H$ be the linear homogeneous polynomial which defines $H\in \sA$ in $\CP^2$. Given a multinet on $\sA$, with parts $(\sA_1,\cdots,\sA_k)$ and
multiplicity vector $m$, set
 $$g_i=\prod_{H\in \sA_i} f_H^{m_H}$$ for $1\leq i\leq k$. 
By definition of multinet, $\{g_i\}_{i=1}^k$ are homogeneous polynomials with the same degree.  
 Then it gives a map $f\colon X(\sA) \to \CP^1$
 by $$f(x)= [g_1(x),g_2(x)].$$
 As shown in \cite{FY}, each $g_i$ is a linear combination of $g_1$ and $g_2$ for any $3\leq i\leq k$, hence 
image of $f$ is exactly  $\Sigma_k=\CP^1-\{ k \text{ points} \}$. In particular, $f$ is a large orbifold fibration since $k\geq 3$ for a multinet.

On the other hand, if $\sB\subset \sA$ is a subarrangement, then the inclusion $X(\sA) \to X(\sB)$ induces an epimorphism $\pi_1(X(\sA))\twoheadrightarrow \pi_1(X(\sB)) $. Moreover, it induces an embedding $\sV^1(\sB) \hookrightarrow \sV^1(\sA)$, see e.g. \cite[Proposition A.1]{Suc14}.
Components of $\sV^1(\sA)$ that are not supported on any proper subarrangement are said to be essential.
 
Following \cite[Theorem 11.4]{Suc14B}, we summarize a series of  works due to Falk, Libgober, Pereira and Yuzvinsky as follows.
\bt \cite{LY,FY,PY,Yuz} \label{thm multinet} Let $\sA$  be  an essential  line arrangement in $\mathbb{CP}^2$.
Then  every  positive-dimensional essential component of $\sV^1(\sA)$ is obtained from $f^*(\sV^1(\pi_1^\orb(\Sigma_k,{\bm \mu})))$ for some orbifold fibration $f\colon X(\sA) \to \Sigma_k$
  and either
\begin{itemize}
\item[(1)] $f$ is small, i.e, $k=2$ and f has at least one multiple fiber, or
\item[(2)] $k=3$ or $4$, and $f$ corresponds to a multinet on $L(\sA)$.
\end{itemize}
\et

\subsection{Pointed multinet}
Theorem \ref{thm structure} shows that $\sV^1(\sA)$ consists of three type components:  the positive-dimensional components passing the origin, the translated positive-dimensional components and the finite union of isolated torsion points. The first one is  well-understood by the multinet structure due to Theorem \ref{thm multinet}.  
The last one is still poorly understood. 
Meanwhile, Denham and Suciu introduced a
combinatorial construction, called pointed multinet \cite{DS}, 
to produce translated 1-dimensional component in $\sV^1(\sA)$.
 As shown by Suciu in \cite{Suc02},  there do exist arrangements $\sA$ with isolated torsion
 points and translated positive-dimensional component in $\sV^1(\sA)$.

Fix a hyperplane $H \in \sA$. The arrangement $\sA'=\sA\setminus\{H\}$ is then called the deletion of $\sA$
with respect to $H$. The pointed multinet structure is a variation of  Definition \ref{def multinet}.
\bd \cite[Definition 5.5]{DS} A pointed multinet on $\sA$ is a multinet structure $((\sA_1,\cdots, \sA_k), m, \mathfrak{X} )$,
together with a distinguished hyperplane $H\in \sA$ for which $m_H>1$, and $m_H \mid  n_x$ for each 
$x\in H\cap \mathfrak{X}$.
\ed
Note that Theorem \ref{thm 34} implies that the pointed multinet structure only exists for $k=3$. Without loss of generality, we assume that $H\in \sA_3$. 
By the multinet structure, we have a large orbifold map $f= [g_1,g_2]\colon X(\sA) \to \Sigma_3$. We can extend this map to $X(\sA')$, denoted by $f'$.  Since $H$ is deleted from $\sA$, the image of $f'$ equals to $\C^*$ and $f'$ has at least one multiple fiber with multiplicity $m_H$. 

\bp \cite[Proposition 5.6]{DS} Suppose that $\sA$ admits a pointed multinet, and $\sA'$ is obtained from $\sA$
by deleting the distinguished hyperplane $H$. Then $\sA'$ supports a small orbifold fibration, and $\sV^1(\sA')$ 
has a component which is a 1-dimensional subtorus, translated by a character of order $m_H$.
\ep 
%\br  \label{rem linear} By this proposition, the translated 1-dimensional subtorus coming from the pointed multiple structure supports a small orbifold fibration
%$f\colon X\to \C^*$ such that at least one multiple fiber of $f$ decomposes as product of linear polynomials.
%\er
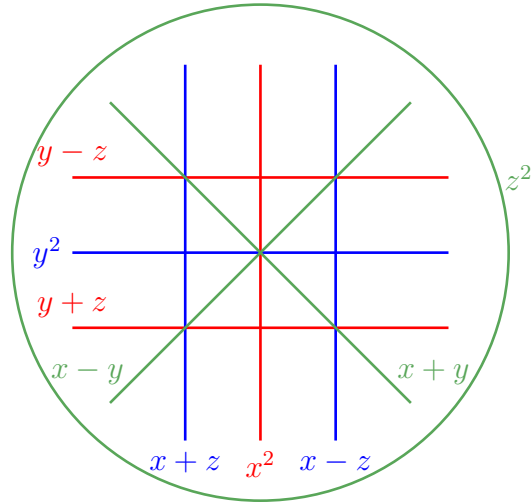
\begin{figure}
\centering
\begin{tikzpicture}
\draw[red,line width=1pt] (0,-2.5) -- (0,2.5);
\node[red,below] at (0,-2.5){$x^2$};
\draw[red,line width=1pt] (-2.5,1) -- (2.5,1);
\node[red,above] at (-2.5,1){$y-z$};
\draw[red,line width=1pt] (-2.5,-1) -- (2.5,-1);
\node[red,above] at (-2.5,-1){$y+z$};

\draw[blue,line width=1pt] (-2.5,0) -- (2.5,0);
\node[blue,left] at (-2.5,0){$y^2$};
\draw[blue,line width=1pt] (-1,-2.5) -- (-1,2.5);
\node[blue,below] at (-1,-2.5){$x+z$};
\draw[blue,line width=1pt] (1,-2.5) -- (1,2.5);
\node[blue,below] at (1,-2.5){$x-z$};

\draw[color=gray!70!green,line width=1pt] (-2,-2) -- (2,2);
\node[color=gray!70!green,above] at (-2.3,-1.9){$x-y$};
\draw[color=gray!70!green,line width=1pt] (-2,2) -- (2,-2);
\node[color=gray!70!green,above] at (2.3,-1.9){$x+y$};
\draw[color=gray!70!green,line width=1pt] (0,0) circle (3.3);
\node[color=gray!70!green,right] at (3.1,1) {$z^2$};
\end{tikzpicture}
\caption{$(3,4)$-multinet on $B_3$-arrangement} \label{Figure 1}
\end{figure}
\bex\cite[Example 3.11]{Suc14}  Let $\sA$ be the reflection arrangement of type $B_3$, defined by the
polynomial $xyz(x^2 -y^2)(x^2 - z^2)(y^2 -z^2)$ with $[x,y,z] \in \CP^2$. Figure \ref{Figure 1} (with partition by colors and multiplicity assigned) gives  a $(3,4)$-multinet on it. The map $$f= [x^2(y^2-z^2), y^2(x^2-z^2)]\colon X(\sA)\to \CP^1$$ has image $\Sigma_3=\CP^1\setminus\{[0,1],[1,0],[1,1]\}.$
It shows that $\sV^1(\sA)$ has a 2-dimensional essential component passing through the origin. 

Let $\sA'$ be the arrangement obtained by deleting the hyperplane $\{z = 0\}$. Then $\sA'$ has the pointed multinet structure with respect to the $(3,4)$-multinet on $\sA$. In fact, the map
 $$f'= [x^2(y^2-z^2), y^2(x^2-z^2)]\colon X(\sA')\to \CP^1$$ has image $\Sigma_2=\CP^1\setminus\{[0,1],[1,0]\}$ and has a multiple fiber at the point $[1,1]$ with multiplicity $2$. 
Hence $\sV^1(\sA')$ has a translated 1-dimensional component. 
\eex 

Suciu asked the following question \cite[Problem 3.12]{Suc14}.
\begin{que} \label{que Suc}
     Do all positive-dimensional translated tori in degree one cohomology jump loci of a line arrangement complement arise in the manner of pointed multinet structure?
\end{que}
%This question can be viewed as a tempt to show that all positive-dimensional component of $\sV^1(X)$ are %combinatorially determined.
 In next section, we show that the double star arrangement introduced by Ishibashi, Sugawara and Yoshinaga in \cite[Example 3.2]{ISY22} 
admits a small orbifold map whose multiple fiber is not formed by lines. This example shows one of the difficult points to prove that the translated positive-dimensional component of $\sV^1(\sA)$ are combinatorially determined:  the multiple fiber of the orbifold map may have irreducible fiber with degree greater than one. 
In particular, it gives a negative answer to Question \ref{que Suc}. 

Note that Theorem \ref{thm multinet} implies that for a line arrangement $\sA$  a translated positive-dimensional component of $\sV^1(\sA)$ may have dimension 1, 2 or 3. The pointed multinet structure can only produce translated 1-dimensional component in $\sV^1(\sA)$. The double star arrangement we will studied later has translated 1-dimensional component in $\sV^1(\sA)$ not coming from the pointed multinet structure. 
It is an interesting question if one can find examples of line arrangement $\sA$ such that $\sV^1(\sA)$ has a translated component with dimension 2 or 3.

\bigskip

\textbf{Acknowledgments.} 
The first author would like to thank Professor Masahiko Yoshinaga for valuable discussion.
Both authors are  supported by National Key Research and Development Project SQ2020YFA070080,  the Project of Stable Support for Youth Team in Basic Research Field, CAS (YSBR-001),  the project ``Analysis and Geometry on Bundles" of Ministry of Science and Technology of the People's Republic of China and  Fundamental Research Funds for the Central Universities.

\section{The double star arrangement}
The double star arrangement consists of 11 lines in $\CP^2$, which is a complexified real arrangement. See Figure \ref{Figure 2} with 5 blue lines, 5 red lines and one green line at the infinity (the big green circle).
Set $\theta=\frac{2}{5}\pi $. We list the defining equations of the ten lines in $\C^2$ as follows: 

$$\begin{aligned}
l_1&=(\sin{2 \theta} - \sin{\theta}) \cdot x + (\cos{2 \theta} - \cos{\theta}) \cdot y - \sin{\theta}\cdot\frac{\cos{\theta}}{\cos{2\theta}}, \\
l_2&=-\sin{\theta} \cdot x + (\cos{\theta} - 1)\cdot y + \sin{\theta}\cdot\frac{\cos{\theta}}{\cos{2\theta}}, \\
l_3&=2(\sin{2 \theta})\cdot x + \sin{\theta}\cdot\frac{\cos{\theta}}{\cos{2\theta}},\\
l_4&=\sin{\theta}\cdot x + (\cos{\theta} - 1)\cdot y - \sin{\theta}\cdot\frac{\cos{\theta}}{\cos{2\theta}}, \\
l_5&=-(\sin{2 \theta} - \sin{\theta})\cdot x + (\cos{2 \theta} - \cos{\theta}) \cdot y + \sin{\theta}\cdot\frac{\cos{\theta}}{\cos{2\theta}}, \\
l_6&=(\sin{2 \theta} - \sin{\theta}) \cdot x + (\cos{2 \theta} - \cos{\theta}) \cdot y - \sin{\theta}, \\
l_7&=-\sin{\theta}\cdot x + (\cos{\theta} - 1)\cdot y + \sin{\theta}, \\
l_8&=2(\sin{2 \theta})\cdot x + \sin{\theta}, \\
l_9&=\sin{\theta}\cdot x + (\cos{\theta} - 1) \cdot y - \sin{\theta}, \\
l_{10} &=-(\sin{2\theta} - \sin{\theta}) \cdot x + (\cos{2 \theta} - \cos{\theta}) \cdot y + \sin{\theta}. 
\end{aligned}$$

\begin{center}
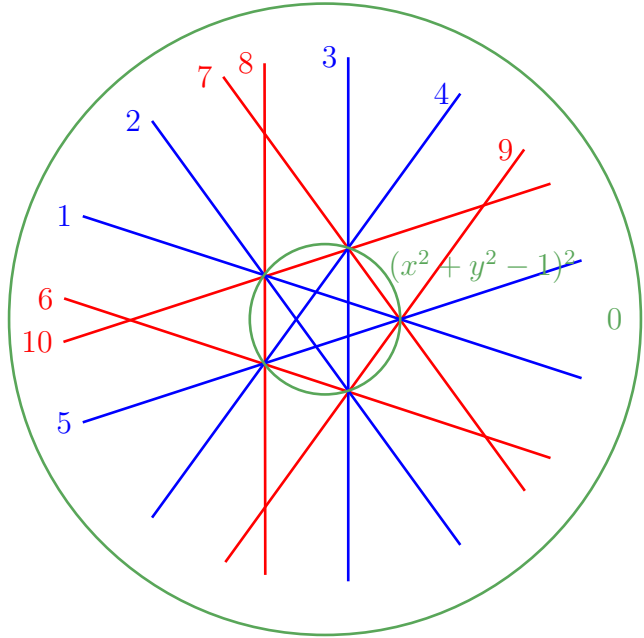

    
\begin{tikzpicture}

\draw[red,line width=1pt] (-3.477,-0.301) -- (2.998,1.804);
\node[red,left] at (-3.477,-0.301){$10$};

\draw[red,line width=1pt] (-3.468,0.2773) -- (2.998,-1.846);
\node[red,left] at (-3.468,0.2773){$6$};

\draw[red,line width=1pt] (-0.8035,3.406) -- (-0.794,-3.401);
\node[red,left] at (-0.8035,3.406){$8$};

\draw[red,line width=1pt] (-1.353, 3.226) -- (2.657,-2.282);
\node[red,left] at (-1.353, 3.226){$7$};

\draw[red,line width=1pt] (-1.325, -3.231) -- (2.648,2.259);
\node[red,left] at (2.648,2.259){$9$};

\draw[blue,line width=1pt] (3.41119, -0.783443) -- (-3.22021, 1.37123);
\node[blue,left] at (-3.22021, 1.37123){$1$};

\draw[blue,line width=1pt] (-3.22021, -1.37123) -- (3.41119, 0.783443);
\node[blue,left] at (-3.22021, -1.37123){$5$};

\draw[blue,line width=1pt] (-2.29921, 2.63887) -- (1.79921, -3.00214);
\node[blue,left] at (-2.29921, 2.63887){$2$};

\draw[blue,line width=1pt] (0.309017, -3.48633) -- (0.309017, 3.48633);
\node[blue,left] at (0.309017, 3.48633){$3$};

\draw[blue,line width=1pt] (-2.29921, -2.63887) -- (1.79921, 3.00214);
\node[blue,left] at (1.79921, 3.00214){$4$};

\draw[color=gray!70!green,line width=1pt] (0,0) circle (4.2);
\draw[color=gray!70!green,line width=1pt] (0,0) circle (1);
\node[color=gray!70!green,right] at (0.7,0.7) {$(x^2+y^2-1)^2$};
\node[color=gray!70!green,right] at (3.6,0) {$0$};

\end{tikzpicture}
\captionof{figure}{The double star arrangement}\label{Figure 2}
\end{center}

\bt \label{thm main}
The double star arrangement is a counterexample to Question \ref{que Suc}.
\et 
\begin{proof} Consider  $X(\sA)$ as the complement of 10 affine lines in $\C^2$.  
Using computations by \textit{Mathematica 12.0} (the codes are listed in the end of this paper), 
we get that the product of the five red lines in $\C^2$ has defining equation 
$$
h_1=c^{-1}\big(\frac{2x^5}{5} -4x^3 y^2 +2xy^4
+\frac{1+\sqrt{5}}{2}x^4 +(1+\sqrt{5})x^2 y^2 + \frac{1+\sqrt{5}}{2}y^4- (2+\sqrt{5})x^2  - (2+\sqrt{5})y^2+\frac{11+5\sqrt{5}}{10} \big)
$$ 
with $c=\frac{32(5+3\sqrt{5})}{125}\sin \frac{2\pi}{5}$ and the product of the five blue lines  in $\C^2$ has defining equation 
$$
h_2=c^{-1}\big(\frac{2x^5}{5} -4x^3 y^2 +2xy^4
+\frac{1-\sqrt{5}}{2}x^4 +(1-\sqrt{5})x^2 y^2 + \frac{1-\sqrt{5}}{2}y^4 - (2-\sqrt{5})x^2  - (2-\sqrt{5})y^2+\frac{11-5\sqrt{5}}{10}\big). 
$$

Consider the map $h= \dfrac{h_1}{h_2} \colon X(\sA) \to \C^*$. It is easy to see that $h$ is surjective and has connected generic fiber.
For any $\lambda\in \C^*\setminus\{1\}$, $h_1-\lambda h_2$ has degree 5. It implies that the map $h$ does not have a multiple fiber over $\lambda \in \C^*\setminus\{1\}$. On the other hand, when $\lambda=1$, 
we have the following equality 
$$
h_1-h_2= \frac{\sqrt{5}}{c} \cdot (x^2+y^2-1)^2.
$$
Let $\overline{h_1}$ and $\overline{h_2}$ be the homogenization of $h_1$ and $h_2$, respectively. 
Then we have \be  \label{multiple}
\overline{h_1}-\overline{h_2}=\frac{\sqrt{5}}{c}z \cdot (x^2+y^2-z^2)^2,
\ee 
where $\{z=0\}$ is the line at infinity.
This shows that the map $h$ has a multiple fiber of multiplicity 2 at the point $1\in \C^*$. 
%\textcolor{red}{In particular, this multiple fiber is not formed by lines. }
So $h$ is a small orbifold fibration with exactly one multiple fiber of multiplicity 2, which gives a translated 1-dimensional component of $\sV^1(\sA)\subseteq \Cha(X(\sA))\cong (\C^*)^{10}$. This translated 1-dimensional component  is of the form $\rho T$, where $\rho=(1,1,1,1,1,-1,-1,-1,-1,-1)$ and $T=(t,t,t,t,t,t^{-1},t^{-1},t^{-1},t^{-1},t^{-1})$. Here we order the affine lines as in Figure \ref{Figure 2}. 

Since the degree 2 polynomial $x^2+y^2-z^2$ appearing in  the equality (\ref{multiple}) can not decompose as a product of two linear  polynomial, this particular translated 1-dimensional component can not come from any pointed multinet structure. 
%It is easy to see that the translated 1-dimensional component induced by $h$ is of the form $\rho T$, where $\rho=(1,1,1,1,1,-1,-1,-1,-1,-1)$ and $T=(t,t,t,t,t,t^{-1},t^{-1},t^{-1},t^{-1},t^{-1})$. Here we order the affine lines as in Figure \ref{Figure 2}
\end{proof}

\br 
 We noticed the double star arrangement due to the following reason. 
For the double star arrangement with lines labeled as in Figure \ref{Figure 2}, 
set $$\omega=\sum_{i=1}^5 \omega_i -\sum_{i=6}^{10} \omega_i,$$ where $\omega_i\in H^1(X(\sA),\Z)$ maps the corresponding meridian to 1 and the meridians associated to other affine lines  to 0.  
Consider the Aomoto complex given by $$0\to H^0(X(\sA),\Z)\overset{\wedge \omega}{\to} H^1(X(\sA),\Z)\overset{\wedge \omega}{\to} H^2(X(\sA),\Z)\to 0$$
According to the computations in \cite[Example 3.2]{ISY22} and the fact that $h$ is small,   the second cohomology of this Aomoto complex must have $2$-torsion. Then \cite[Proposition 1.12]{LL21} shows that the map $ h$ may have a multiple fiber with multiplicity dividing by $2$.
\er

One can use the following codes in \textit{Mathematica 12.0} to check that the equations of $h_1$ and $h_2$ listed in the proof of Theorem \ref{thm main}  give Figure \ref{Figure 2}. 
\begin{lstlisting}
Print["h1=",h1=2/5 *x^5-4x^3*y^2+2x*y^4+(1+Sqrt[5])/2 *x^4+(1+Sqrt[5])x^2*y^2+(1+Sqrt[5])/2 *y^4-(2+Sqrt[5])x^2-(2+Sqrt[5])y^2+(11+5Sqrt[5])/10];
Print["h2=",h2=2/5 *x^5-4x^3*y^2+2x*y^4+(1-Sqrt[5])/2 *x^4+(1-Sqrt[5])x^2*y^2+(1-Sqrt[5])/2 *y^4-(2-Sqrt[5])x^2-(2-Sqrt[5])y^2+(11-5Sqrt[5])/10];
ContourPlot[{h1==0,h2==0},{x,-3.5,3.5},{y,-3.5,3.5},PlotLegends->{"Expressions",{"h1","h2"}}]
\end{lstlisting}

\bibliographystyle{amsalpha}

\end{document}